\documentclass[letterpaper,10pt,conference]{ieeeconf}
\IEEEoverridecommandlockouts  
\usepackage{wrapfig}
\usepackage{float}
\usepackage[framemethod=TikZ]{mdframed}
\usepackage{relsize}

\usepackage{amssymb,amsmath,verbatim,bm,color,times,mathtools,multicol}
\usepackage{graphicx}
\usepackage{psfrag}
\usepackage{hyperref} 
\usepackage{cleveref}
\usepackage{siunitx}
\usepackage{tabularx}
\usepackage{multirow,array,booktabs}
\usepackage[justification=centering]{caption}
\usepackage{subcaption}
\usepackage{wrapfig}
\usepackage{textcomp}
\usepackage{algorithm,algorithmic}
\usepackage[font=small,textfont=it]{caption}
\usepackage{stfloats}
\usepackage{import}
\usepackage{xifthen}
\usepackage{pdfpages}
\usepackage[left=48pt, right=48pt, bottom=76pt, top=50pt]{geometry}

\usepackage[style=ieee]{biblatex}

\usepackage[colorinlistoftodos]{todonotes}

\newcommand{\diag}{\mathrm{diag}}

\newcommand\ddfrac[2]{\frac{\displaystyle #1}{\displaystyle #2}}

\begin{document}

\title{\Huge \bf Flatness-based MPC using B-splines transcription with application to \\ a Pusher-Slider System}
% \author{Thomas Neve\textsuperscript{1,2,*}, Tom Lefebvre\textsuperscript{1,2} and Guillaume Crevecoeur\textsuperscript{1,2}
\author{Thomas Neve\textsuperscript{1,2,*}, Tom Lefebvre\textsuperscript{1,2}, Sander De Witte\textsuperscript{1,2} and Guillaume Crevecoeur\textsuperscript{1,2}
    \thanks{
        \textsuperscript{1}D2Lab research group, Department of Electromechanical, Systems and Metal Engineering, Ghent University, Tech Lane Ghent Science Park 913, B-9052 Zwijnaarde, Belgium}
    \thanks{\textsuperscript{2}Core Lab MIRO, Flanders Make.}
    \thanks{\textsuperscript{*}Corresponding author: \texttt{thomas.neve@ugent.be}.
    }
    \vspace*{-25pt}
}

\maketitle

\begin{abstract}
    This work discusses the use of model predictive control (MPC) for the manipulation of a pusher-slider system.
    In particular we leverage the differential flatness of the pusher-slider in combination with a B-splines transcription to address the computational demand that is typically associated to real-time implementation of an MPC controller.
    We demonstrate the flatness based B-spline MPC controller in simulation and compare it to a standard MPC implementation approach using direct multiple shooting. We evaluate the computational advantage of the flatness based MPC empirically and document computational acceleration up to 65\%.
\end{abstract}

% \listoftodos

\section{Introduction}
\noindent
% The field of robotic manipulation has advanced beyond static pick-and-place tasks, allowing robots to perform real-time sensor-based manipulation tasks in human environments \cite{stuber2020let}.
When we manipulate objects in a daily setting our actions reach further than just simple grasping and releasing.
We perform an array of different tasks, such as pushing, pulling, sliding, and lifting.
Therefore, we argue that pushing is a motion primitive of practical significance for robotic manipulation as well.
Pushing extends the normal capabilities resulting in a wide variety of applications, such as the manipulation of hard to grasp objects due to for example their weight or shape, or positioning tasks \cite{yu2016more}.
Moreover, we contend that pushing an object allows to simplify the design of the end-effector without significantly compromising the maneuverability of the object.
The pusher-slider system, discussed in this work, represents a basic non-prehensile manipulation task where the goal is to control the motion of the slider through the pusher.
The system consists of a sliding object (the slider) and a single contact point (the pusher).

To achieve optimal control of a pusher-slider system, a precise model of the slider's behavior under pushing is required, as well as a robust control strategy to address model inaccuracies and external disturbances.
Model Predictive Control (MPC) has proven to be an effective control strategy in complex systems due to its ability to anticipate future behavior and act optimally in changing environments \cite{hogan}.
MPC optimizes dynamic behavior over a prediction horizon at each sampling period while adhering to a set of constraints.
Common approaches include direct methods, which transforms the optimal control problem into a numerical optimization problem \cite{Diehl2009}.
Several techniques exist here to transcribe the control problem, typically using polynomial or piecewise constant parameterization.
A common direct approach is Direct Multiple Shooting (DMS) which considers both the state and control as optimization variables.
Subsequent states in the prediction horizon are tied together using equality constraints between each control interval \cite{Albersmeyer2009}.

%to simpler systems or requiring significant effort from a control engineer to develop an efficient MPC controller.
% There are several approaches to address this challenge, such as approximating the system dynamics to reduce complexity or transcribing the Optimal Control Problem (OCP) in an efficient manner for faster solution times with existing solvers.
% In this work we focus on the first approach, approximating the system dynamics and using the resulting differentially flat properties of its approximated model.

Solving this optimization problem in real-time can be computationally demanding, limiting its practical application.
In this work we attempt to adress this challenge by exploitation of the differentially flat properties of the dynamics model.
Differential flatness is a property of a class of dynamical systems characterized by the ability to define all states and controls of the system as a set of specific differential variables and their derivatives.
This property can be particularly useful for both solving trajectory planning \cite{stoican2015flat} and tracking problems \cite{nguyen2017, Greeff2018, wang2019}.
In \cite{HELLING202014686} full optimal flatness based MPC has also been used for the control of underactuated surface vessels.
Most work however does not solve the full optimal control control problem each MPC iteration.
It is still unclear how the differential flatness property is best used in an MPC implementation and what the impact is on computational requirement.
% This work tackles the following question: \textit{How do we leverage the differential flatness property in an MPC implementation to decrease computational demand?}

We present an MPC controller that leverages the differentially flat properties of the introduced kinematic pusher-slider model and employs a B-spline to parameterize the trajectory.
% We attempt to answer this by employing B-splines in the transcription of the optimal control problem to optimization problem.
This significantly reduces the complexity of the trajectory optimization problem due to less optimization variables and constraints compared to the DMS approach.
The B-spline is used to represent the trajectory of the flat coordinates of the pusher-slider system, after which the other variables can be deduced from the flat coordinates.
However, the nonlinear flat expressions used to compute the original states and controls from the flat trajectory also add some nonlinearity to the objective and constraints.
It is unclear how such a flatness based MPC would compare to standard transcription methods such as DMS.
Our contribution consists of the application of a flatness based MPC strategy to a pusher-slider system, and comparing its performance in terms of computational requirements to the DMS approach.

\section{Pusher-Slider Model}
The pusher-slider system can be modelled as a quasi-static model where no accelerations occur. This model is physically valid in the following cases.
\begin{enumerate}
    \item The motions of the slider are slow enough that inertial forces are negligible compared to friction forces.
    \item The friction forces at the contact point between pusher and slider are negligible with respect to the friction forces between the slider and the ground.
\end{enumerate}

We start by describing the kinematic model and continue by discussing the differentially flat properties of the resulting model and the accompanying expressions of the pusher-slider system.

\subsection{Kinematics}
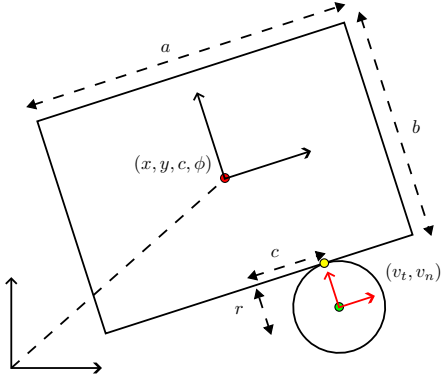
\begin{figure}
    \centering
    \def\svgwidth{\columnwidth}
    \scalebox{0.7}{%% Creator: Inkscape 1.2.1 (9c6d41e410, 2022-07-14), www.inkscape.org
%% Accompanies image file 'schematic.pdf' (pdf, eps, ps)
\begingroup%
  \makeatletter%
  \providecommand\color[2][]{%
    \errmessage{(Inkscape) Color is used for the text in Inkscape, but the package 'color.sty' is not loaded}%
    \renewcommand\color[2][]{}%
  }%
  \providecommand\transparent[1]{%
    \errmessage{(Inkscape) Transparency is used (non-zero) for the text in Inkscape, but the package 'transparent.sty' is not loaded}%
    \renewcommand\transparent[1]{}%
  }%
  \providecommand\rotatebox[2]{#2}%
  \newcommand*\fsize{\dimexpr\f@size pt\relax}%
  \newcommand*\lineheight[1]{\fontsize{\fsize}{#1\fsize}\selectfont}%
  \ifx\svgwidth\undefined%
    \setlength{\unitlength}{285.10235495bp}%
    \ifx\svgscale\undefined%
      \relax%
    \else%
      \setlength{\unitlength}{\unitlength * \real{\svgscale}}%
    \fi%
  \else%
    \setlength{\unitlength}{\svgwidth}%
  \fi%
  \global\let\svgwidth\undefined%
  \global\let\svgscale\undefined%
  \makeatother%
  \begin{picture}(1,0.79731523)%
    \lineheight{1}%
    \setlength\tabcolsep{0pt}%
    \put(0,0){\includegraphics[width=\unitlength,page=1]{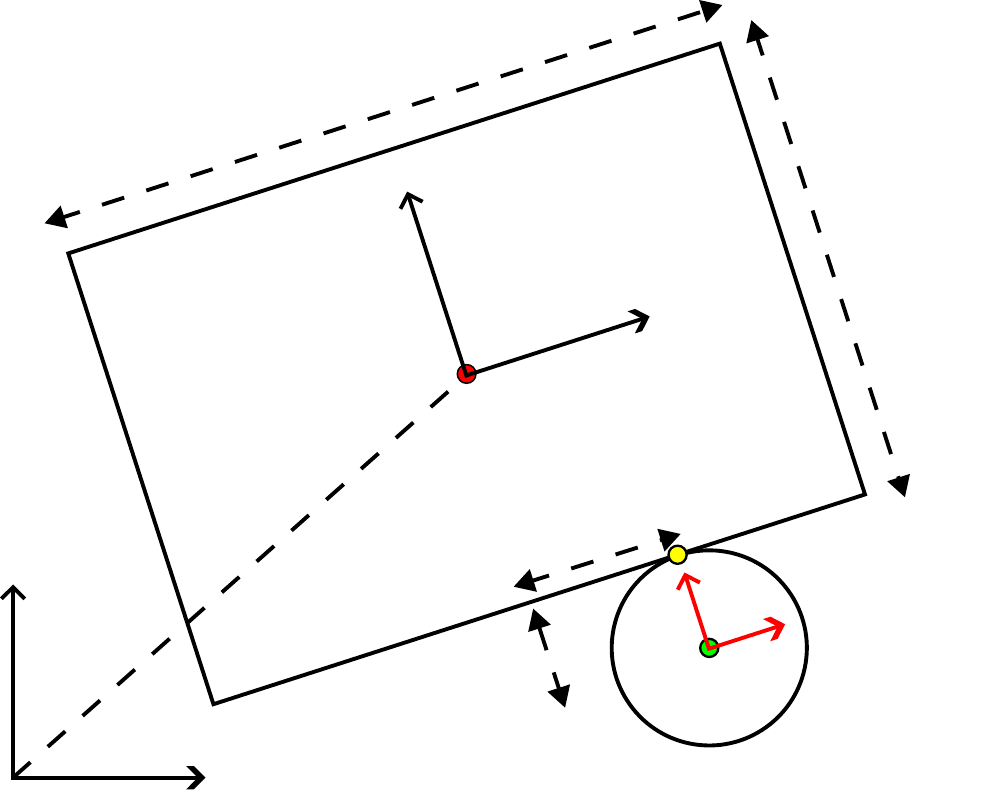}}%
    \put(0.57852133,0.24966839){\color[rgb]{0,0,0}\makebox(0,0)[t]{\lineheight{1.25}\smash{\begin{tabular}[t]{c}$c$\end{tabular}}}}%
    \put(0.50186447,0.12831623){\color[rgb]{0,0,0}\makebox(0,0)[t]{\lineheight{1.25}\smash{\begin{tabular}[t]{c}$r$\end{tabular}}}}%
    \put(0.34403355,0.68943231){\color[rgb]{0,0,0}\makebox(0,0)[t]{\lineheight{1.25}\smash{\begin{tabular}[t]{c}$a$\end{tabular}}}}%
    \put(0.36039608,0.43917444){\color[rgb]{0,0,0}\makebox(0,0)[t]{\lineheight{1.25}\smash{\begin{tabular}[t]{c}$(x,y,c,\phi)$\end{tabular}}}}%
    \put(0.87611366,0.20706322){\color[rgb]{0,0,0}\makebox(0,0)[t]{\lineheight{1.25}\smash{\begin{tabular}[t]{c}$(v_t,v_n)$\end{tabular}}}}%
    \put(0.88166054,0.51678662){\color[rgb]{0,0,0}\makebox(0,0)[t]{\lineheight{1.25}\smash{\begin{tabular}[t]{c}$b$\end{tabular}}}}%
  \end{picture}%
\endgroup%
}

    \caption{State and input values pusher-slider system.}
    \label{fig:pusher_slider_model}
\end{figure}

The system consists of a sliding object (the slider) and a single contact point (the pusher). We define the state of the system as
\begin{equation*}
    \textbf{x}^\top=
    \begin{pmatrix}
        x & y & c & \phi
    \end{pmatrix}
\end{equation*}
and the input as
\begin{equation*}
    \textbf{u}^\top = \begin{pmatrix}
        v_t & v_n
    \end{pmatrix}
\end{equation*}
with variables shown in figure \ref{fig:pusher_slider_model}. Here $(x, y, \phi)$ represents the planar configuration of the slider in the global frame of reference and c the relative position of the pusher contact with the slider.
The inputs $v_t$ and $v_n$ denote, respectively, the speed tangent to and normal to the surface against which is pushed.
Some additional geometric parameters are defined as well: with $a$ and $b$ denoting the length of the slider's side, and $r$ denoting the radius of the pusher.

By assuming quasi-static interaction and frictionless contact between the pusher and the slider, the model becomes entirely kinematic.
The former means physically that we assume that in any configuration of the slider the ground friction is large enough to keep the slider stationary if not acted upon.
The following equations emerge
\begin{align*}
    \dot{x}    & = -\frac{\beta^{2}}{\beta^{2}+c^{2}}v_{n}\sin(\phi)                \\
    \dot{y}    & = \frac{\beta^{2}}{\beta^{2}+c^{2}}v_{n}\cos(\phi)                 \\
    \dot{c}    & = v_{t} - \left(\frac{b}{2}+r\right)\frac{c}{\beta^{2}+c^{2}}v_{n} \\
    \dot{\phi} & = \frac{c}{\beta^{2}+c^{2}}v_{n}
\end{align*}
where $\beta^2$ is a factor that depends on the geometry of the slider.
In the case of a rectangular geometry $\beta^2 = \ddfrac{1}{12}(a^2+b^2)$.
For the complete derivation we refer to the work by Lefebvre et al.\cite{lefebvre}.

\subsection{Differential flatness}
Differential flatness or just flatness is a generalisation of the notion of inverse dynamics for underactuated nonlinear systems, hence $n_x > n_u$ where $n_x$ and $n_u$ denote the state $\textbf{x} \in \mathbb{R}^{n}$ and control $\textbf{u} \in \mathbb{R}^{m}$ dimensionality.
To circumvent the over-defined inverse dynamics problem, a set of flat coordinates $\textbf{y} \in \mathbb{R}^{n_u}$ is selected with the same dimensionality as the control.
The flat coordinates themselves may not have physical meaning.
\[\textbf{y} = y(\textbf{x}, \textbf{u}, \dot{\textbf{u}}, \ddot{\textbf{u}},...,\textbf{u}^{(p)})\]
Then the system state and the control are expressed as a function of the flat coordinates and its derivatives.
\[\textbf{x} = x( \textbf{y}, \dot{\textbf{y}}, \ddot{\textbf{y}},...,\textbf{y}^{(q)})\]
\[\textbf{u} = u( \textbf{y}, \dot{\textbf{y}}, \ddot{\textbf{y}},...,\textbf{y}^{(q)})\]

It can be shown that the equations of motion of the pusher-slider are flat \cite{lefebvre}.
A set of flat coordinates with the same number of coordinates as the input $n_u = 2$ is defined.
For the pusher-slider system the flat coordinates are taken as the Cartesian coordinates.
\[
    \textbf{p}^\top = \begin{pmatrix}
        x & y
    \end{pmatrix}
\]
The flat expressions for the control input and auxiliary states can be derived as a function of these coordinates and their derivatives.
\begin{align}
    \label{eq:flat_expressions}
    \begin{split}
        c    & = \beta^2 \ddfrac{\dot{x}\ddot{y}-\ddot{x}\dot{y}}{\sqrt{\dot{x}^2+\dot{y}^2}^3}                                                                                                                                        \\
        \phi & = - \arctan \ddfrac{\dot{x}}{\dot{y}}                                                                                                                                                                                   \\
        v_t  & = \left(1+\beta^2\ddfrac{(\ddot{x}\dot{y} - \dot{x}\ddot{y})^2}{(\dot{x}^2+\dot{y}^2)^3}\right) \sqrt{\dot{x}^2+\dot{y}^2}                                                                                                         \\
        v_n  & =  \beta^2 \ddfrac{\dot{x}\dddot{y}-\dddot{x}\dot{y}}{\sqrt{\dot{x}^2+\dot{y}^2}^3} + 3 \beta^2 \ddfrac{(\ddot{x}\dot{y} - \dot{x}\ddot{y})(\dot{x}\ddot{xy}+\dot{y}\ddot{y})}{\sqrt{\dot{x}^2+\dot{y}^2}^5} + \dotsm \\
        & \quad \left(\ddfrac{b}{2}+r\right)\ddfrac{\dot{x}\ddot{y}-\ddot{x}\dot{y}}{\dot{x}^2+\dot{y}^2}
    \end{split}
\end{align}
% where $\vectorstyle{u}$ is expressed in local coordinates.
% To obtain an expression in global coordinates one can use the following transformation {\color{red} Is this relevant here, maybe leave out?}
% \[
%     \vectorstyle{u}^g = \matrixstyle{R}_{\phi} \vectorstyle{u}^l = \ddfrac{1}{\sqrt{\dot{x}^2+\dot{y}^2}}\begin{bmatrix}
%         \dot{y}  & \dot{x} \\
%         -\dot{x} & \dot{y}
%     \end{bmatrix} \vectorstyle{u}^l
% \]

\section{Model Predictive Control}
Model predictive control (MPC) is a control strategy that optimizes over a prediction horizon to steer a dynamic system in an optimal way whilst satisfying a set of constraints.
% The major challenge in real-time MPC is balancing these requirements with the computation requirement.

Consider the following optimal control problem (OCP)
\begin{align}
    \label{eq:continuous_OCP}
    \begin{split}
        \min_{\textbf{u}(t)} l_T(\textbf{x}(T)) + \int\nolimits_t^T l(\textbf{x}(t), \textbf{u}(t)) dt
        \\[5pt]
        \text{s.t.}
        \left\{
        \begin{array}{ll}
            \dot{\textbf{x}}(t) = f(\textbf{x}(t),\textbf{u}(t)) \\
            \textbf{x}(0) = \textbf{x}_{\text{start}}            \\
            \textbf{x}(T) = \textbf{x}_{\text{goal}}             \\
            \textbf{h}(\textbf{x}(t), \textbf{u}(t)) \leq \textbf{0}
        \end{array}
        \right.
    \end{split}
\end{align}
with $\textbf{x}(0)$, $\textbf{x}(T)$ and $\textbf{h}$ the initial, final state and path constraints respectively.
Each MPC iteration, the OCP is solved to acquire the optimal control signal $\textbf{u}(t)$.
This control signal could be applied in open-loop to the system to perform in an optimal manner.
However, due to model inaccuracies or external disturbances the system could quickly diverge from the expected behavior.
To account for this, the actual state following the applied control is used as the new initial state, and the optimization process is repeated. Effectively closing the loop.
Typically, the control signal $\textbf{u}(t)$ is considered over the same time horizon $T$ on each iteration.
Should one want to converge to a specific goal state within the given time horizon, it is also possible to reduce this horizon on each iteration.
Thus modifying the OCP formulation on each step.

MPC could be used for both trajectory tracking and generation. The former uses an MPC controller to track a predefined trajectory \cite{Greeff2018,wang2019}.
In this work however we focus on trajectory generation, generating a new trajectory at each time step.
Ergo we reconsider the full optimal control problem every time step.

% A common MPC approach is a direct method.
% Here the open-loop control problem is transformed into a numerical optimization problem, also referred to as trajectory optimization, by discretizing the dynamic model.
% Common approaches include direct methods, which transform the optimal control problem into a numerical optimization problem.
The common approach to solving the optimal control problem in MPC is through the use of direct methods, which translate the control problem into a numerical optimization problem.
Within the class of direct methods several approaches exist to transcribe the OCP into a trajectory optimization problem.
A popular approach for longer time horizons is direct multiple shooting (DMS), which considers both the system state and control as optimization variable.
Another approach is the use of polynomial parameterization, e.g. B-splines \cite{HELLING202014686}, of the trajectory.
This approach can be particularly usefull in combination with a differentially flat system.
We continue this section with a more detailed discussion on both approaches.
These techniques will then be applied and compared onto the pusher-slider system in terms of computational performance to highlight their respective advantages and disadvantages.

\subsection{Direct Multiple Shooting}
Direct multiple shooting (DMS) is a strategy to transcribe a trajectory optimization problem into a numerical optimization problem. The prediction horizon of the problem is split, using a discretized dynamics model, into several equidistant shooting nodes.
The multiple shooting transcription of the control problem (\ref{eq:continuous_OCP}) results in the following nonlinear program (NLP)
\begin{align}
    \label{eq:traj_opt}
    \begin{split}
        \min_{\textbf{u}_{\text{0:N-1}},\textbf{x}_{\text{0:N}}}\sum_{i=0}^{N-1}\,l(\textbf{x}_{i},\textbf{u}_{i})
        \\
        \text{s.t.}
        \left\{
        \begin{array}{ll}
            \textbf{x}_{i+1} = f(\textbf{x}_i,\textbf{u}_i)                          \\
            \textbf{x}_{0} = \textbf{x}_{\text{start}}                               \\
            \textbf{x}_{N} = \textbf{x}_{\text{goal}}                                \\
            \textbf{x}_{\text{min}} \leq \textbf{x}_{i} \leq \textbf{x}_{\text{max}} \\
            \textbf{u}_{\text{min}} \leq \textbf{u}_{i} \leq \textbf{u}_{\text{max}}
        \end{array}
        \right.
    \end{split}
\end{align}
with $N$ the number of discretization steps, $l$ the cost function and $\textbf{f}$ the discrete dynamics of the pusher-slider over one shooting interval. % where each node consists of a state at the start of the interval and a zero-order hold control during the node interval.
Note the addition of the state trajectory $x_t$ as an optimization variable.
This component is key to the mutliple shooting approach.
Continuity in the state trajectory between shooting nodes is enforced through a set of continuity constraints using $f$.

The addition of the states $x_t$ to the optimization might seem counterproductive compared to a single shooting approach where only the controls $u_t$ are considered as optimization variables.
But compared to single shooting, mutliple shooting is more flexible in initalizing the problem, and has improved convergence properties \cite{Albersmeyer2009, quirynen2015}.  % second citation calls first

\subsection{Flatness based MPC}
The continuity constraint of (\ref{eq:traj_opt}) is implicitly fulfilled in the flat parameterization (\ref{eq:flat_expressions}) of the system.
Put differently, the entire class of feasible state-action trajectories is encoded by the class of smooth flat paths.
Therefore, the differential flatness properties of the pusher-slider could be leveraged to transcribe the control problem without the need for a set of constraints to enforce continuity.
To parametize the state trajectory, B-splines can be used to represent the flat trajectory without considering the controls themselves as optimization variables directly.
From this parameterized trajectory the full state and controls can still be inferred using (\ref{eq:flat_expressions}).

\subsubsection{B-splines}
In this section we give an introduction to B-splines.
For an extensive theory we refer to the related literature \cite{de1978practical,Schoenberg1946ContributionsTT}.

B-splines, first introduced in \cite{Schoenberg1946ContributionsTT}, consist of a union of local curve segments which are each active on a specific interval. % A spline function consists of several polynomial segments defined over subintervals which are joined together using continuity constraints.
Its segmented nature allows for very efficient tailoring to desired local changes \cite{Ilkiv}.
This property also makes it particularly useful for trajectory optimization where one might desire a local change without affecting global behaviour, which in the end will result in sparse Jacobian and Hessian structures of the nonlinear program.

Consider the expansion
\begin{equation}
    \label{eq:bspline_def}
    S(\tau) = \sum_{i=0}^{n} N_{i,j} p_i, \quad  \tau \in [0, T]
\end{equation}
where the spline $S(\tau)$ is defined over the closed interval $[0, T]$ which is subdivided into $m$ sub intervals with $m=k+n+1$. The vector $\textbf{p}=[p_0, \dots, p_n]$ constitutes the set of control points which acts as a set of weights on the basis functions $N_{i,j}$ where $k$ is the order of the B-spline, $n$ is the number of control points and $m$ defining the number of knots.
The basis functions can be determined according to the following recursion formula \cite{DEBOOR197250}.
\begin{equation}
    \label{eq:recursion_formula}
    \begin{aligned}
        N_{i,0}(\tau) = \begin{cases}
                            1 & \tau_i \leq \tau \leq \tau_{i+1} \\
                            0 & \text{otherwise}
                        \end{cases},
        \\
        N_{i,j}(\tau) = \frac{\tau - \tau_i}{\tau_{i+j}-\tau_i} N_{i,j-1}(\tau) \\
        + \frac{\tau_{i+j+1}-\tau}{\tau_{i+j+1}-\tau_{i+1}} N_{i+1, j-1}(\tau)
    \end{aligned}
\end{equation}
It can be seen here that the closed interval $[0, T]$ is divided uniformly through a knot vector.
\begin{equation}
    \boldsymbol{\tau} = [\tau_0, \tau_1, \dots, \tau_m]
\end{equation}
At the knot points, the polynomials are joined and connected in a continuous manner.
From the recursion (\ref{eq:recursion_formula}) it is clear that each basis function $N_{i,j}$, with weight $p_i$, is only active on a subset of the interval $[0, T]$.
More specifically, it is nonzero on the interval $[\tau_i, \tau_{i+p+1})$, resulting in the B-spline not being defined at the start and the ending of the interval $[0, T]$.
Choosing the knot vector with duplicate knots at the ends as
\begin{equation}
    \boldsymbol{\tau} = [\underbrace{0 \dots 0}_\text{k} \underbrace{0 \dots T}_\text{internal knots} \underbrace{T \dots T}_\text{k}],
\end{equation}
alleviates this and also clamps the start and endpoint of the spline to the two end control points, $S(0)=p_0$ and $S(T)=p_n$.
This is also often referred to as a clamped uniform B-spline. % uniform referring to the distribution of the knot vector

The derivative of a B-spline is simply a function of B-splines of a lower degree.
Since a B-spline function of order $k$ consists of polynomials of order $n-1$, it's derivatives are continuous up to the derivative of degree $n-2$.
The flat expressions require a derivative of the flat coordinates up to order $q$.
Thus the order of the B-spline used to represent the flat trajectory should be at least of order $q+2$.
This is usually chosen as the minimum value to avoid numerical issues \cite{HELLING202014686}.

\subsubsection{B-spline transcription}
We can now represent our OCP as a numerical optimization problem using B-spline transcription.
Here we consider the control points $\textbf{p}_{\text{0:n}}$ as optimization variables that bend the continuous B-spline curve according to some objective and constraints.
With an equidistant knot vector of length $N+1$, resulting in a uniform B-spline, we can evaluate the B-spline function at a set of collocation $N+1$ points.
Using the flat expressions (\ref{eq:flat_expressions}) we can express the objective and constraints at the collocation points, resulting in an optimization problem without the need for continuity constraints.
\begin{align}
    \label{eq:B-spline_ocp}
    \begin{split}
        &\min_{\textbf{p}_{\text{0:n}}}\sum_{i=0}^{N}\,
        l(
        \textbf{x}_{i}(\textbf{p}_{\text{0:n}}) ,
        \textbf{u}_{i}(\textbf{p}_{\text{0:n}})
        )
        \\
        &\text{s.t.}
        \left\{
        \begin{array}{ll}
            \textbf{x}_{0}(\textbf{p}_{\text{0:n}}) = \textbf{x}_{\text{start}}                               \\
            \textbf{x}_{N}(\textbf{p}_{\text{0:n}}) = \textbf{x}_{\text{goal}}                                \\
            \textbf{x}_{\text{min}} \leq \textbf{x}_{i}(\textbf{p}_{\text{0:n}}) \leq \textbf{x}_{\text{max}} \\
            \textbf{u}_{\text{min}} \leq \textbf{u}_{i}(\textbf{p}_{\text{0:n}}) \leq \textbf{u}_{\text{max}}
        \end{array}
        \right.
    \end{split}
\end{align}

\section{Simulation Experiments}
% We tested several cases where the slider starts and moves from the initial position, $[x, y]=[0, 0]$, and orientation, $\phi=0$, to a certain goal pos and orientation.
We will now verify the computational characteristics of the proposed transcription methods in function of real-time appliction of the corresponding MPC approach.
Therefore we define several validation cases.
Due to the arbitrariness of the global frame of reference we can consider initial position $[x, y]=[0, 0]$ and orientation $\phi=0$ without loss of generality.
Each case the MPC controller is used to control the pusher to manoeuvre the slider according to some objective towards the goal location and orientation.
Each MPC iteration an OCP is solved at discrete time points in time with a fixed time step $\Delta t = \frac{T}{N+1}$.
Two different transcription techniques, DMS and B-spline transcription, are used to build the MPC controller. For each transcription technique we consider both a fixed and shrinking horizon version, resulting in four different implementations.
The controllers are implemented in CasADi \cite{Andersson2019} and solved using IPOPT \cite{Wachter2006}.
In the simulation, noise was added to each control, $u_{noise} \sim \mathcal{N}(0, \Sigma)$, with $\Sigma=\diag(0.2, 0.2)$.

More concretely, the OCP takes the following generic form for both the spline and DMS transcription.
\begin{align}
    \label{eq:OCP}
    \begin{split}
        J = \min_{\textbf{u}_{0:N-1}}&\sum_{i=0}^{N-1}\,\textbf{u}_{i}^T R \textbf{u}_{i} + \textbf{x}_{i}^T Q \textbf{x}_{i} + \textbf{x}_{N}^T P \textbf{x}_{N}
        \\
        &\text{s.t.}
        \left\{
        \begin{array}{ll}
            \textbf{x}_{0} = \textbf{x}_{\text{start}} \\
            -0.4 \leq c \leq 0.4
        \end{array}
        \right.
    \end{split}
\end{align}

Next to the objective desribed above, we add a regularization term to the B-spline transcription case, $J_{\text{tot}}=J+w_{\text{reg}} J_{\text{reg}}$.
\begin{equation}
    J_{\text{reg}} = \sum_{i=0}^{n-1} \|\textbf{p}_{i+1}-\textbf{p}_{i} \|^2
\end{equation}
This term was added with the aim of maintaining a uniform distance between the control points.
We have found empirically that this can improves numerical stability significantly.

Initial results here were promising for the most part.
However, once the slider got close to the goal, MPC performance starts deteriorating in the form of increasing computation time and numerical instability.
This phenomenon ocurred for both the B-spline and DMS approach but was more outspoken for the flatness based MPC which starts to fail slightly sooner.
For this reason, the MPC controller is used up untill a switching point where the final number of steps are handled by a tracking controller.
The tracking controller uses a DMS formulation to track the final solution from the MPC generation.

% We use two seperate criteria here for the fixed and shrinking horizon.

We start by describing our methodology of an MPC with typical fixed preview horizon.

\subsection{Fixed horizon}
Each MPC step an OCP with preview horizon $T=1$ and $N=20$ discretization steps is solved.
The first control from the solution, $u_0$, is applied to the system and aformentioned OCP is resolved with the new state, position and orientation, of the pusher-slider.

% The MPC controller is used from the initial state on for a number of steps up untill a switching point.
% Once the slider starts getting close to the goal, $\|(x, y) - (x_{\text{goal}}, y_{\text{goal}}) \| = 2.5$, the controller is switched to a tracking controller.

% \subsubsection{B-spline}
% The solution resulting from the B-spline transcription is a contnuous trajectory in $\tau$.
% In order to sample the control we need to assume some relation between time $t$ and $\tau$ used to parameterize the spline.
% Surprisingly, the states of the pusher-slider $[x, y, \phi, c]$ are independent of the chosen time dependency of the path.
% In contrast to the expressions for the states, the controls, $[v_t, v_n]$, depend on the parameterization $\tau(t)$.
% If however the relation is assumed linear, the change in the expressions is straightforward \cite{lefebvre} with $\dot{\tau}=\frac{1}{T}$.
% \begin{equation}
%     \begin{aligned}
%         v_t = \dot{\tau} v_t(\tau) \\
%         v_n = \dot{\tau} v_n(\tau)
%     \end{aligned}
%     \label{eq:control_computation}
% \end{equation}
% The first control can then be computed at $\tau=0$, and applied to the simulation.

% In our experiments we use a B-spline of order $k=5$ with number of knots $m=12$ resulting in $n=6$ control points.

% \subsubsection{DMS}
% Extracting the first control, $u_0$, is straightforward by taking the control on the first shooting node.
\subsubsection{Implementation details}
The MPC controller is used up untill the slider starts getting close to the goal, $\|(x, y) - (x_{\text{goal}}, y_{\text{goal}}) \| = 2.5$. Afterwards the tracking controller is actived.

The solution resulting from the B-spline transcription is a continous trajectory in $\tau$.
The required control can be easily computed with the expressions (\ref{eq:flat_expressions}) using the flat coordinates and its derivatives evaluated at $\tau=0$.

In our experiments we use a B-spline of order $k=5$ with number of knots $m=12$ resulting in $n=6$ control points.

\subsubsection{Discussion}
\begin{figure}
    \centering
    \includegraphics[width=\linewidth]{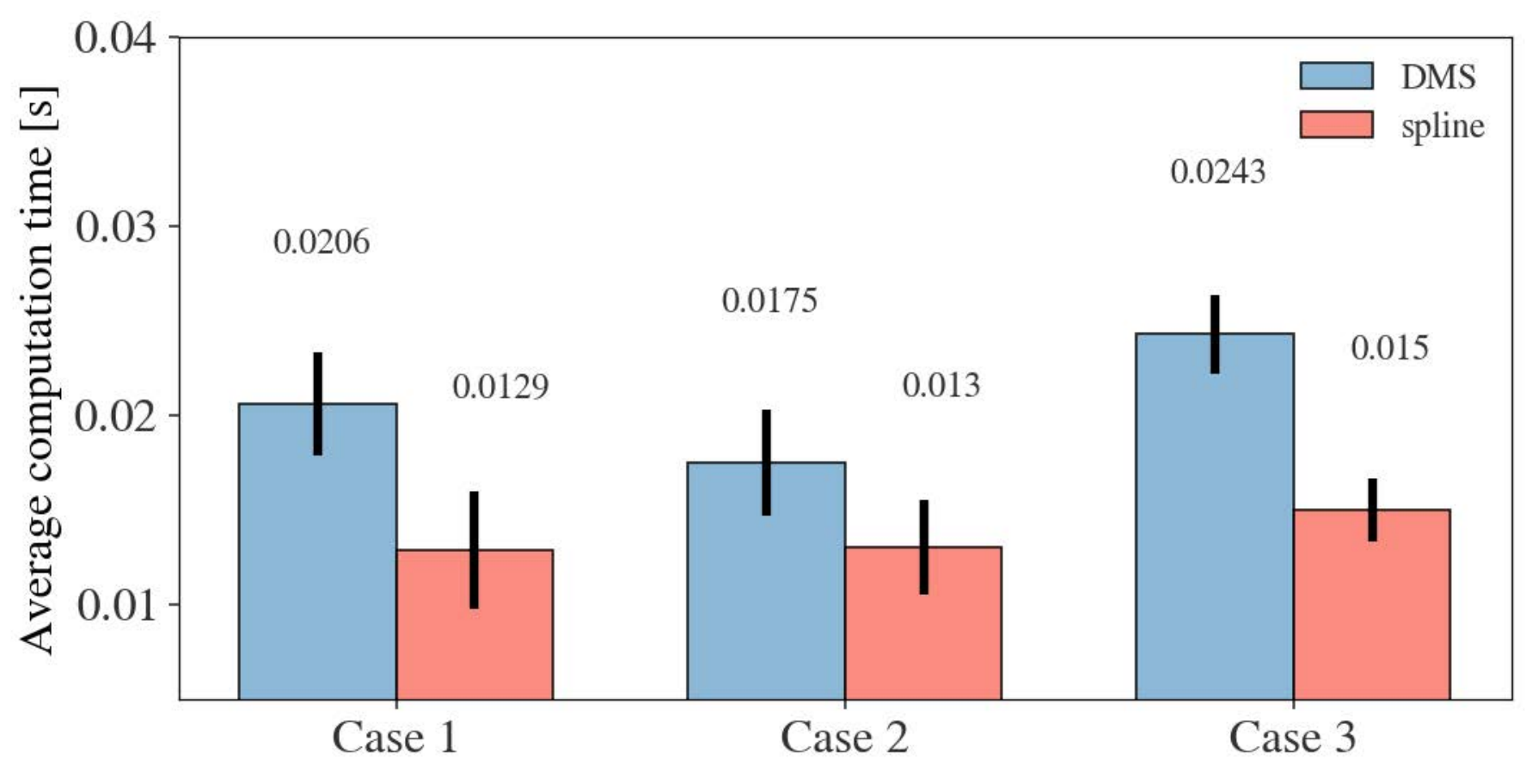}
    \caption{Computation time of an MPC controller with fixed preview horizon for all cases.}
    \label{fig: computation scen 1}
\end{figure}

The results for three test cases are shown in figure \ref{fig: case 1, fixed}, \ref{fig: case 2, fixed} and \ref{fig: case 3, fixed} with the asociated computation times indicated in figure \ref{fig: computation scen 1}.
The MPC phase is indicated in green and tracking phase of the control is indicated in red.
The indicated computation times are only given for the MPC phase.

In all cases, the B-spline transcription with flat trajectories outperforms the DMS method in terms of computation time.
The main reason for this is the substantially lower amount of optimization variables in the B-spline transcription.
Only $n$ control points are used, resulting in $n \cdot n_u$ optimization variables.
In contrast, the DMS approach requires a control and state vector for each of the $N$ shooting nodes, resulting in $N \cdot n_x \cdot n_u + n_x$ variables.

% \begin{figure}
%     \centering
%     \includegraphics[width=\linewidth]{1_scen1.pdf}
%     \caption{case 1: fixed horizon, settings: $R=\diag(2, 2)$, $Q=\diag(1, 1, 10, 0)$, $P=\diag(500, 500, 2000, 0)$, and $w_{\text{reg}}=2$}
%     \label{fig: case 1, fixed}
% \end{figure}
% \begin{figure}
%     \centering
%     \includegraphics[width=\linewidth]{3_scen1.pdf}
%     \caption{case 2: fixed horizon, settings: $R=\diag(0.5, 2)$, $Q=\diag(1, 1, 0, 0)$, $P=\diag(500, 500, 2000, 0)$ and $w_{\text{reg}}=10$}
%     \label{fig: case 2, fixed}
% \end{figure}
% \begin{figure}
%     \centering
%     \includegraphics[width=\linewidth]{4_scen1.pdf}
%     \caption{case 3: fixed horizon, settings: $R=\diag(2, 2)$, $Q=\diag(0.5, 0.5, 0, 0)$, $P=\diag(500, 500, 2000, 0)$ and $w_{\text{reg}}=5$}
%     \label{fig: case 3, fixed}
% \end{figure}

\begin{figure}
    \centering
    \begin{subfigure}{0.49\linewidth}
        \includegraphics[width=\linewidth]{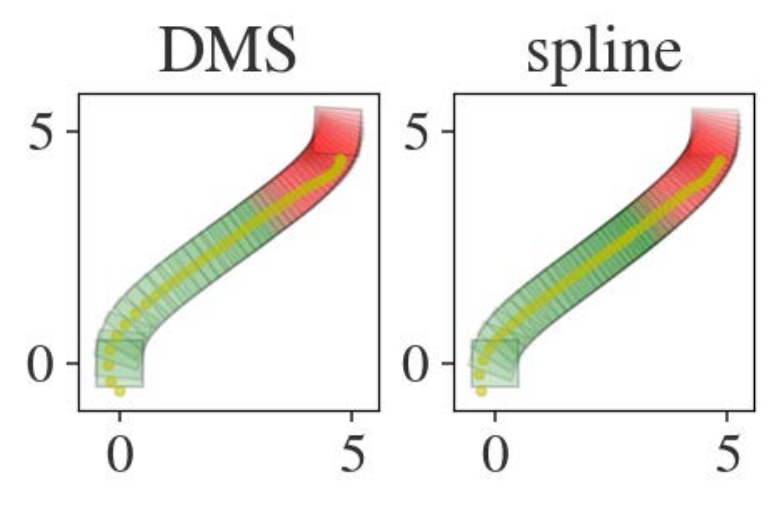}
        % \caption{case 1, $R=\diag(2, 2)$, $Q=\diag(1, 1, 10, 0)$, $P=\diag(500, 500, 2000, 0)$, and $w_{\text{reg}}=2$}
        \caption{case 1}
        \label{fig: case 1, fixed}
    \end{subfigure}
    \hfill
    \begin{subfigure}{0.49\linewidth}
        \includegraphics[width=\linewidth]{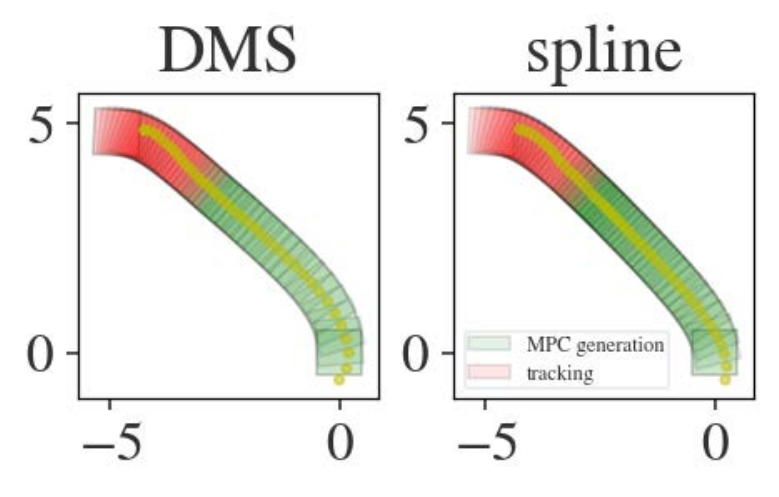}
        % \caption{case 2, $R=\diag(2, 2)$, $Q=\diag(0.5, 0.5, 0, 0)$, $P=\diag(500, 500, 2000, 0)$ and $w_{\text{reg}}=5$}
        \caption{case 2}
        \label{fig: case 2, fixed}
    \end{subfigure}
    \begin{subfigure}{\linewidth}
        \includegraphics[width=\linewidth]{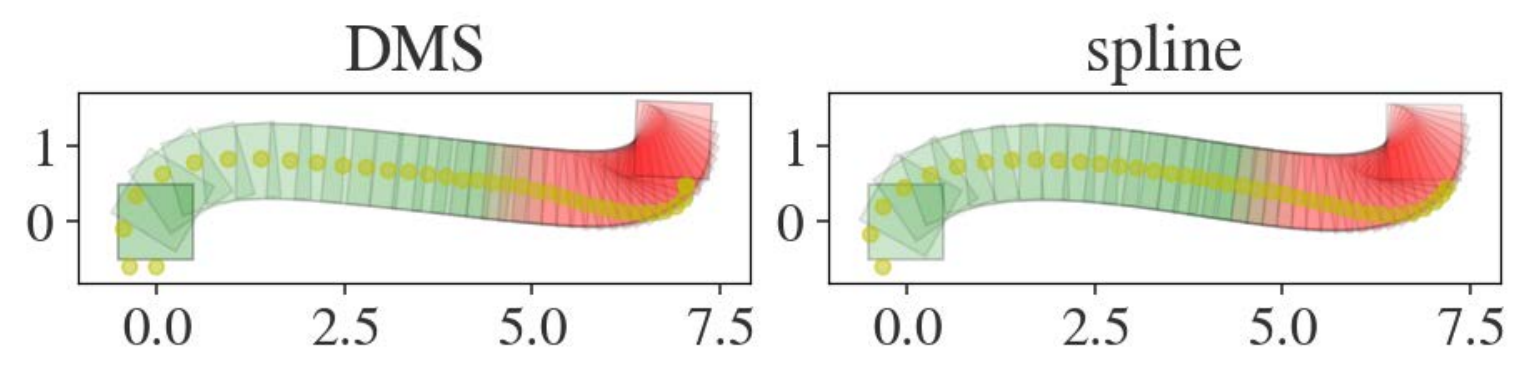}
        % \caption{case 3, $R=\diag(0.5, 2)$, $Q=\diag(1, 1, 0, 0)$, $P=\diag(500, 500, 2000, 0)$ and $w_{\text{reg}}=10$}
        \caption{case 3}
        \label{fig: case 3, fixed}
    \end{subfigure}
    \caption{All cases with a fixed horizon MPC. The settings are reported in table \ref{table:settings_fixed_horizon}.}
    % Case 1: $R=\diag(2, 2)$, $Q=\diag(1, 1, 10, 0)$, $P=\diag(500, 500, 2000, 0)$, $w_{\text{reg}}=2$, case 2: $R=\diag(2, 2)$, $Q=\diag(0.5, 0.5, 0, 0)$, $P=\diag(500, 500, 2000, 0)$, $w_{\text{reg}}=5$ and case 3:$R=\diag(0.5, 2)$, $Q=\diag(1, 1, 0, 0)$, $P=\diag(500, 500, 2000, 0)$, $w_{\text{reg}}=10$
\end{figure}
\begin{table}
    \centering
    \begin{tabular}{|c | c | c | c | c |}
        \hline
        settings & $\diag$(R)   & $\diag$(Q)         & $\diag$(P)            & $w_{\text{reg}}$ \\
        \hline
        case 1   & $(2, 2)$     & $(1, 1, 10, 0)$    & $(500, 500, 2000, 0)$ & 2                \\
        case 2   & $(2, 2)$     & $(0.5, 0.5, 0, 0)$ & $(500, 500, 2000, 0)$ & 5                \\
        case 3   & $(0.5, 0.5)$ & $(1, 1, 0, 0)$     & $(500, 500, 2000, 0)$ & 10               \\
        \hline
    \end{tabular}
    \caption{Settings used for the fixed horizon MPC. For the cost matrices, the diagonal components are reported.}
    \label{table:settings_fixed_horizon}
\end{table}

\subsection{Shrinking horizon}
Instead of keeping the horizon fixed, each subsequent step the horizon is reduced with the goal of converging to the goal within the intended number of steps $N=20$ and time horizon $T=1$.
In other words, after each control step, the preview horizon $T$ is reduced by $\Delta t$.

\subsubsection{Implementation details}
In contrast to the fixed horizon formulation, we also add a terminal constraint to the problem, $\textbf{x}_N=\textbf{x}_{\text{goal}}$.
This change in the formulation also allows us to drop the Mayer term $P$ and $Q$ from the lagrangian term in (\ref{eq:OCP}).
Relying solely on the terminal constraint and reducing horizon to converge to the goal.
The tracking controller is also actived after $N_{\text{MPC}}=15$ steps.

Since the B-spline transcription results in a continuous flat geometric trajectory, redefining the time parameterization is very practical.
One only needs to replace the time horizon $T$ in the B-spline parameterization (\ref{eq:bspline_def}) with the new reduced horizon.
The number of optimization variables also remains identical throughout the manoeuvre.

% \subsubsection{DMS}
For the DMS approach, a reducing horizon is conceptually relatively simple but can take a bit more effort in implementation.
Each control step the number of steps in the preview horizon $N$ is reduced by 1 in order to retain the same time step $\Delta t$ under a reducing time horizon.
This also results in a reducing number of variables and constraints on each subsequent iteration.
Here we opted to rebuild the underlying NLP on each iteration.

\subsubsection{Discussion}
\begin{figure}
    \centering
    \includegraphics[width=\linewidth]{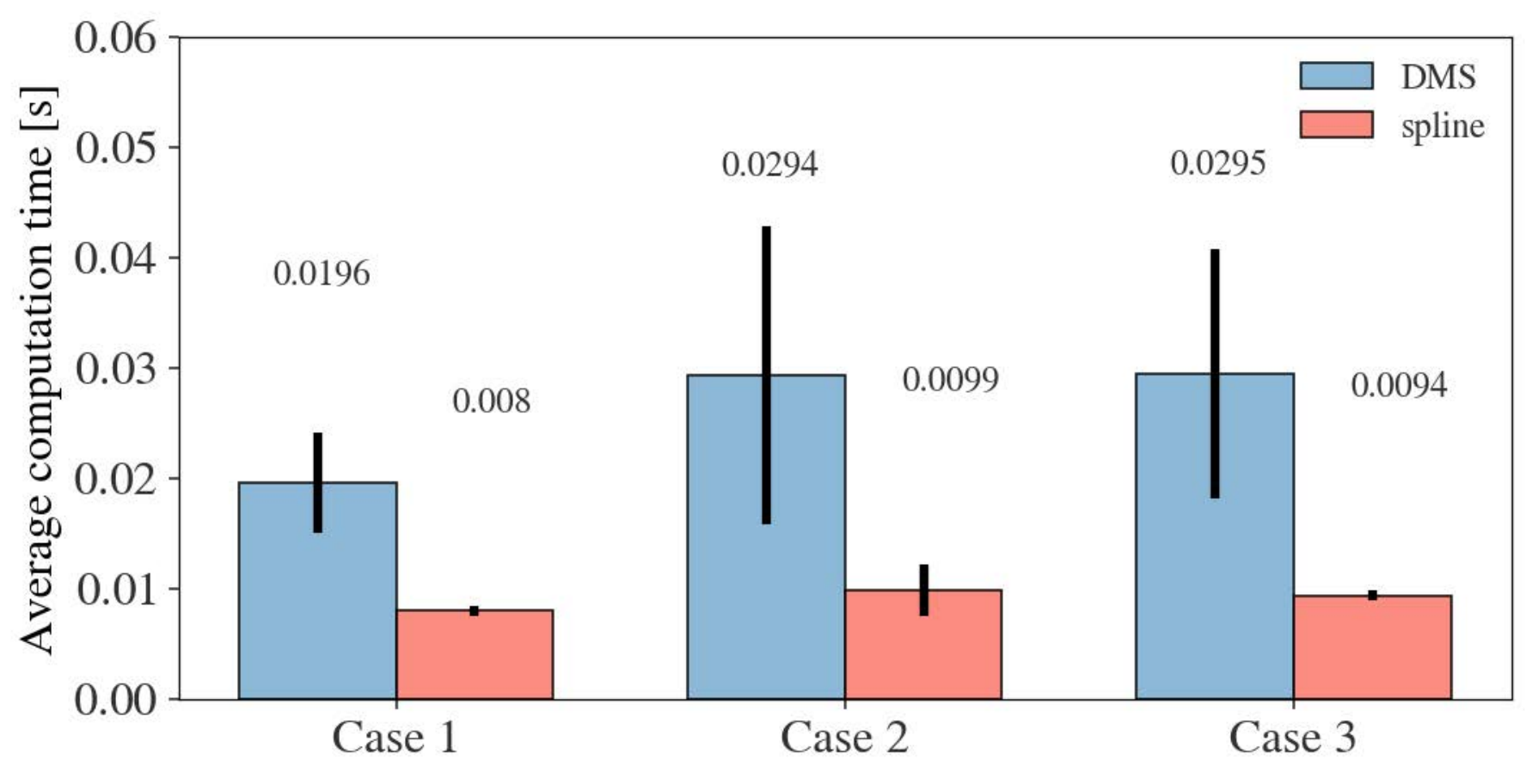}
    \caption{Computation time of an MPC controller with shrinking horizon for all cases.}
    \label{fig: computation scen 2}
\end{figure}

The results for three test cases are shown in figure \ref{fig: case 1, shrinking}, \ref{fig: case 2, shrinking} and \ref{fig: case 3, shrinking} with the asociated computation times indicated in figure \ref{fig: computation scen 2}.
Similar to the fixed horizon scenario, the B-spline method outperforms the DMS method in computation times.
Notice the higher variance in the computation time of the DMS approach.
This is unsurprising as the number of optimization variables reduces when the number of horizon steps is reduced, resulting in a reduced computational demand upon nearing the goal.
% However, in the DMS approach there is a clear deacreasing trend in computational demand upon nearing the goal.
% \begin{figure}
%     \centering
%     \includegraphics[width=\linewidth]{1_scen2.pdf}
%     \caption{case 1: shrinking horizon, settings: $R=\diag(0.1, 0.1)$, $w_{\text{reg}}=1$}
%     \label{fig: case 1, shrinking}
% \end{figure}
% \begin{figure}
%     \centering
%     \includegraphics[width=\linewidth]{3_scen2.pdf}
%     \caption{case 2: shrinking horizon, settings: $R=\diag(0.4, 0.4)$, $w_{\text{reg}}=10$}
%     \label{fig: case 2, shrinking}
% \end{figure}
% \begin{figure}
%     \centering
%     \includegraphics[width=\linewidth]{4_scen2.pdf}
%     \caption{case 3: shrinking horizon, settings: $R=\diag(0.5, 0.5)$, $w_{\text{reg}}=1$}
%     \label{fig: case 3, shrinking}
% \end{figure}
\begin{figure}
    \centering
    \begin{subfigure}{0.49\linewidth}
        \includegraphics[width=\linewidth]{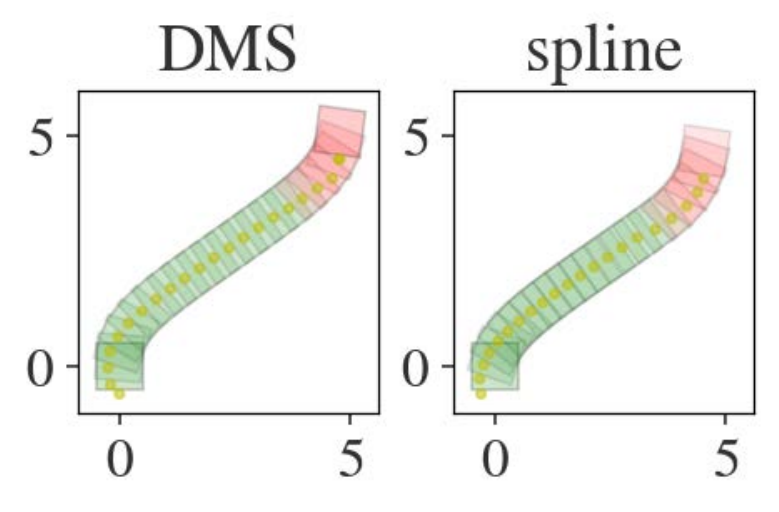}
        % \caption{case 1, $R=\diag(0.1, 0.1)$, $w_{\text{reg}}=1$}
        \caption{case 1}
        \label{fig: case 1, shrinking}
    \end{subfigure}
    \hfill
    \begin{subfigure}{0.49\linewidth}
        \includegraphics[width=\linewidth]{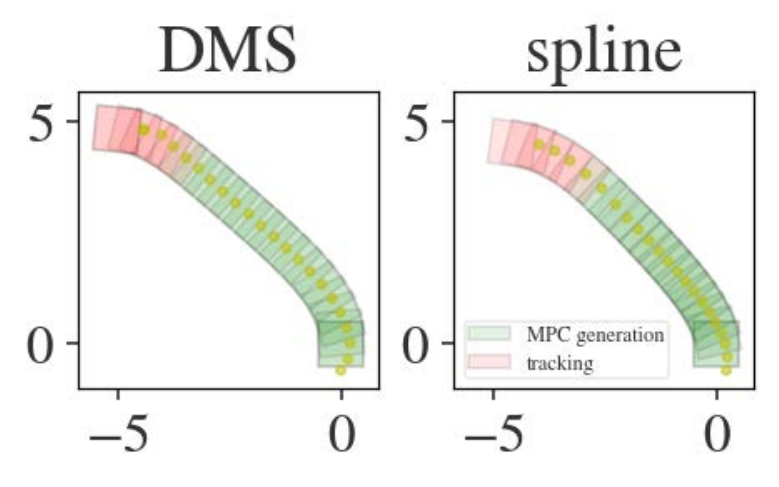}
        % \caption{case 2, $R=\diag(0.5, 0.5)$, $w_{\text{reg}}=1$}
        \caption{case 2}
        \label{fig: case 2, shrinking}
    \end{subfigure}
    \begin{subfigure}{\linewidth}
        \includegraphics[width=\linewidth]{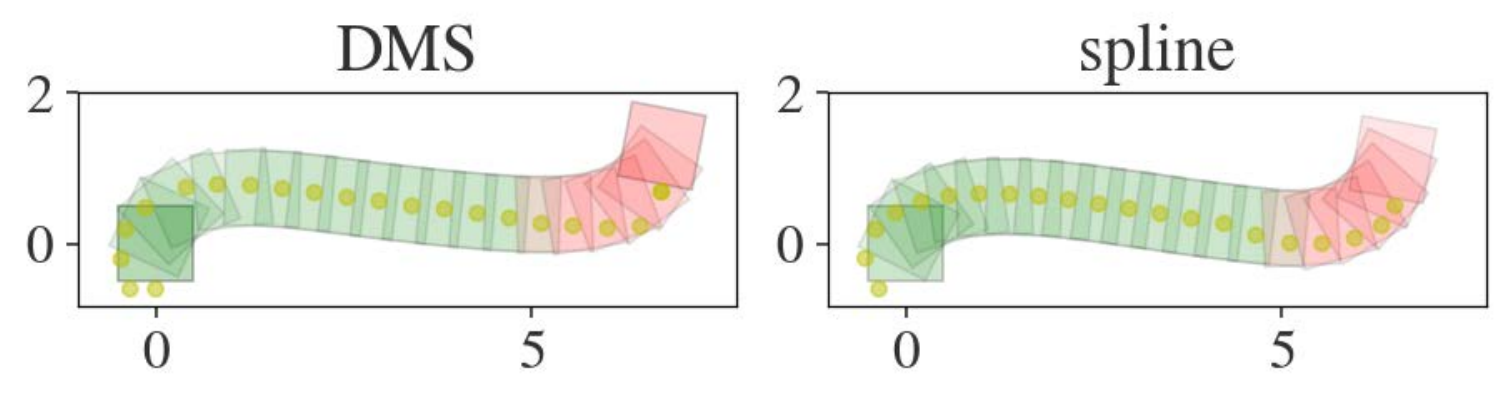}
        % \caption{case 3, $R=\diag(0.4, 0.4)$, $w_{\text{reg}}=10$}
        \caption{case 3}
        \label{fig: case 3, shrinking}
    \end{subfigure}
    \caption{All cases with a shrinking horizon horizon MPC. The settings are reported in table \ref{table:settings_shrinking_horizon}.}
    % Case 1: $R=\diag(0.1, 0.1)$, $w_{\text{reg}}=1$, case 2: $R=\diag(0.5, 0.5)$, $w_{\text{reg}}=1$ and case 3: $R=\diag(0.4, 0.4)$, $w_{\text{reg}}=10$.
\end{figure}

\begin{table}
    \centering
    \begin{tabular}{|c | c | c |}
        \hline
        settings & $\diag$(R)   & $w_{\text{reg}}$ \\
        \hline
        case 1   & $(0.1, 0.1)$ & 1                \\
        case 2   & $(0.5, 0.5)$ & 1                \\
        case 3   & $(0.4, 0.4)$ & 10               \\
        \hline
    \end{tabular}
    \caption{Settings used for the shrinking horizon MPC. For the cost matrices, the diagonal components are reported.}
    \label{table:settings_shrinking_horizon}
\end{table}

\section{Conclusion}
In this work we demonstrated the use of a flatness based MPC using a B-splines transcription to parameterize a flat trajectory on a pusher-slider system.
We compared this MPC methodology to a typical MPC controller using DMS and evaluated how both methodologies compare on a several cases where a slider manoeuvre towards a goal is performed.
We did this for both a fixed and shrinking horizon implementation.

Results show that there is benefit in exploiting the differential flatness of the pusher-slider in combination with B-splines.
One of it's advantages includes a reduced amount of optimization variables required to transcribe the optimal control problem to an optimization problem.
Furthermore, the continuous parameterization of the flat trajectory offers a more elegant way to implement a shrinking horizon MPC.

\section*{Acknowledgment}
This work was supported by the Flanders Make projects DIRAC and the ``Onderzoeksprogramma Artificiële Intelligentie (AI) Vlaanderen'' programme.

% \bibliographystyle{IEEEtran}
% % \bibliographystyle{unsrt}
% \bibliography{bibliography}
\printbibliography

\end{document}